\documentclass[conference]{IEEEtran}
\ifCLASSINFOpdf
   \usepackage[pdftex]{graphicx}
\else
   \usepackage[dvips]{graphicx}
\fi
\usepackage[cmex10]{amsmath}
\usepackage{amsfonts,amssymb,amsthm}
\usepackage{xspace}
\usepackage[linktocpage=true,colorlinks=true,linkcolor=blue,citecolor=blue,urlcolor=blue]{hyperref}
\usepackage{cleveref}
\begin{document}
\title{
An Augmented Lagrangian Method on GPU for Security-Constrained AC Optimal Power Flow
}


\author{\IEEEauthorblockN{François Pacaud\IEEEauthorrefmark{1},
Armin Nurkanović\IEEEauthorrefmark{2},
Anton Pozharskiy \IEEEauthorrefmark{2},
Alexis Montoison\IEEEauthorrefmark{3} and
Sungho Shin\IEEEauthorrefmark{4}}
\IEEEauthorblockA{\IEEEauthorrefmark{1} Centre Automatique et Systèmes, Mines Paris-PSL, Paris, France}
\IEEEauthorblockA{\IEEEauthorrefmark{2} Department of Microsystems Engineering (IMTEK), University of Freiburg, Germany}
\IEEEauthorblockA{\IEEEauthorrefmark{3} Mathematics and Computer Science Division, Argonne National Laboratory, Lemont, IL, USA}
\IEEEauthorblockA{\IEEEauthorrefmark{4} Massachusetts Institute of Technology, Cambridge, Massachusetts, USA }
}
\thanks{
AE  is with TU Dortmund University, Dortmund, Germany (alexander.engelmann@ieee.org). SS is with Argonne National Laboratory, Lemont, IL, USA (sshin@anl.gov). FP is with Centre Mathématiques et Systèmes, Mines Paris-PSL, Paris, France (francois.pacaud@minesparis.psl.eu). VZ is with University of Wisconsin-Madison, Madison, WI, USA and Argonne National Laboratory, Lemont, IL, USA (victor.zavala@wisc.edu).}

\maketitle

\begin{abstract}
  We present a new algorithm for solving large-scale security-constrained optimal power flow in polar form (AC-SCOPF). The method builds on Nonlinearly Constrained augmented Lagrangian (NCL), an augmented Lagrangian method in which the subproblems are solved using an interior-point method. NCL has two key advantages for large-scale SC-OPF. First, NCL handles difficult problems such as infeasible ones or models with complementarity constraints. Second, the augmented Lagrangian term naturally regularizes the Newton linear systems within the interior-point method, enabling to solve the Newton systems with a pivoting-free factorization that can be efficiently parallelized on GPUs. We assess the performance of our implementation, called MadNCL, on large-scale corrective AC-SCOPFs, with complementarity constraints modeling the corrective actions. Numerical results show that MadNCL can solve AC-SCOPF with 500 buses and 256 contingencies fully on the GPU in less than 3 minutes, whereas Knitro takes more than 3 hours to find an equivalent solution.
\end{abstract}

\begin{IEEEkeywords}
  AC-SCOPF;
  Augmented Lagrangian method;
  Contingency screening;
  GPU acceleration;
  Nonlinear programming.
\end{IEEEkeywords}


\section{Introduction}
\subsection{Motivation}
In transmission networks, the optimal dispatch is usually computed by solving
a security-constrained optimal power flow (SCOPF). The dispatch
minimizes a given criterion (costs or network losses) while considering
physical constraints (power flow, line flow limits) and the generators' capacities.
Furthermore, the dispatch should remain feasible
for a set of contingency scenarios corresponding to the loss of a line or a generator in the network ($N-1$ security criterion).
We refer to \cite{stott2005security,frank2016introduction} for comprehensive descriptions of the SC-OPF problem.

The SCOPF is usually formulated as a large-scale linear program called the DC-SCOPF, whose size
grows linearly with the number of contingencies~\cite{alsac2002further}.
This comes at the price of linearizing the nonlinear physical constraints, incurring the solution's accuracy~\cite{coffrin2012approximating}.
However, solving the AC-SCOPF with the original nonlinear formulation remains an open challenge,
and was the main motivation behind the recent GO competition~\cite{aravena2023recent}.
The issue with the nonlinear formulation is two-fold. First, the AC-SCOPF writes as a huge-scale
nonlinear program, whose size is beyond the capacity of state-of-the-art nonlinear optimization solvers
like Ipopt or Knitro~\cite{wachter2006implementation,waltz2006interior,capitanescu2011state}.
Second, the complementarity constraints in the AC-SCOPF models are numerically
difficult to handle, and a tailored algorithm must be employed.
These complementarity constraints arise from the set of logical conditions
modeling the droop control and the PV/PQ switches constraining the recourse in the post-contingency state.
Complementarity constraints have been adopted early
on in power system computation, both for power flow~\cite{zhao2008pv,sundaresh2014modified,murray2015robust} and for optimal power flow~\cite{rosehart2005optimal}.
However, the AC-SCOPF reformulated with complementarity constraints yields a degenerate nonlinear program
coming with its own issues.

\subsection{Augmented Lagrangian method}

There exists numerous algorithms to solve mathematical programs with complementarity constraints (MPCCs) \cite{nurkanovic2024}.
For instance, the teams in the GO-competition has handled the complementarity constraints
using active-set methods~\cite{curtis2023decomposition}, homotopy reformulations~\cite{gholami2023admm,petra2023surrogate} or an $\ell_1$-exact penalty reformulation~\cite{leyffer2006interior}.
The large-scale nature of AC-SCOPFs pushes algorithms to their boundary, as the number
of complementarity constraints grows linearly with the number of contingencies.
To address this challenge, we propose to solve the AC-SCOPF with an augmented Lagrangian
method based on the Nonlinearly Constrained augmented Lagrangian (NCL)~\cite{ma2017stabilized,montoison2025madncl}. The augmented Lagrangian comes with three benefits for SCOPF problems:
(1) it can quickly detect local infeasibility (a common situation for AC-SCOPF)~\cite{chiche2016augmented},
(2) it is robust and can handle complementarity constraints~\cite{izmailov2012global},
(3) it has a structure favorable for GPU acceleration, as the Newton systems we obtain
in the algorithm can be solved efficiently in parallel without numerical pivoting~\cite{montoison2025madncl}.
MadNCL is a recent implementation of Algorithm NCL built on top of MadNLP, and supports
the solution of large-scale problems on the GPU using the linear solver NVIDIA cuDSS.
As such, MadNCL~\cite{montoison2025madncl} can be considered as an extension of our previous work investigating
the solution of large-scale optimal power flow (OPF) on the GPU with MadNLP~\cite{shin2024accelerating}.


\subsection{Scope and contributions}
In this work, we analyze the performance of MadNCL~\cite{montoison2025madncl}
on large-scale AC-SCOPF instances formulated as MPCCs.
To the best of our knowledge, this is the first time an augmented Lagrangian-based algorithm is used to solve the corrective SCOPF formulated with complementarity constraints. We present numerical
experiments showing that MadNCL can accommodate well the complementarity constraints in the AC-SCOPF and
is effective at detecting infeasible problems --- a useful feature for contingency screening.
We detail how to deport the computation on the GPU for faster solution time.
On the GPU, MadNCL evaluates the derivatives with ExaModels~\cite{shin2024accelerating} and solves the Newton systems with NVIDIA cuDSS.
We compare MadNCL with Artelys Knitro~\cite{waltz2006interior}, a state-of-the-art optimization
solver that supports the solution of MPCCs~\cite{leyffer2006interior}.

\section{Model}
In this section, we detail how to formulate the AC-SCOPF with complementarity constraints.
We adapt the formulation used in the GO competition~\cite{aravena2023recent}.
We denote by $K$ the number of contingencies, and refer to the base case by the index $k=0$.
We suppose the system has $n_b$ buses, $n_\ell$ lines and $n_g$ generators.
For $k=0, \cdots, K$, the voltage magnitudes and angles at buses are denoted by
$(v_b^k, \theta_b^k) \in \mathbb{R}^{n_b} \times \mathbb{R}^{n_b}$, and the active and reactive power generations by $(p_g^k, q_g^k) \in \mathbb{R}^{n_g} \times \mathbb{R}^{n_g}$.
For $(a, b) \in \mathbb{R}^2$, we say that $a$ complements $b$ if $a \geq 0$, $b \geq 0$
and $a b = 0$. The complementarity constraint is denoted by $0 \leq a \perp b \geq 0$.

\subsection{Recourse constraints}
In SCOPF, the variation of the power production $p_g^k \in \mathbb{R}^{n_g}$ in contingency $k$
should reflect the behavior of the \emph{automatic generation control system}
(AGC, also known as droop control): the active power is used to regulate the
frequency in the post-contingency state. Given a participation factor
encoded as a vector $\alpha_g \in \mathbb{R}^{n_g}$, the power generation
in contingency $k$ is given by
\begin{equation}
  \label{eq:agc}
  p_g^k = \min\Big(\max\big( p_g^0 + \alpha_g \Delta^k , \; \underline{p}_g \big), \; \overline{p}_g \Big) \; ,
\end{equation}
where $\Delta^k \in \mathbb{R}$ is a variable encoding the power adjustment in contingency $k$
and $\underline{p}_g, \overline{p}_g$ are two vectors encoding the lower and upper bounds on the
active power generation.
The non-smooth $\min$ and $\max$ operations clip $p_g^0 + \alpha_g \Delta^k$ to the bounds, and rewrite equivalently as
a set of complementarity constraints~\cite{baumrucker2008mpec}: for non-negative $\pi_{g,+}^k \geq 0$
and $\pi_{g,-}^k \geq 0$, equation~\eqref{eq:agc} is equivalent to,
for all $k=1, \cdots, K$,
\begin{equation}
  \label{eq:droopmpec}
  \begin{aligned}
    & \pi_{g,+}^k - \pi_{g,-}^k = p_{g}^k - (p_g^0 + \alpha_g \Delta) \; , \\
    & 0 \leq \pi_{g,-}^k \perp \overline{p}_g - p_g^k \geq 0 \; , \\
    & 0 \leq \pi_{g,+}^k \perp p_g^k - \underline{p}_g \geq 0 \; .
  \end{aligned}
\end{equation}
Similarly, the \emph{voltage control} keeps the voltage
magnitudes at the PV buses at their nominal values $v_{m,b}^k = v_{m,b}^0$ by injecting
or absorbing reactive power. However, if the reactive power production $q_g^k$ reaches
its lower $\underline{q}_g$ or upper limit $\overline{q}_g$, the bus is converted to a PQ bus and the voltage is allowed
to vary. The PV/PQ switches are modeled with a second set of complementarity constraints writing,
for $\nu_{b-}^k \geq 0$ and $\nu_{b+}^k \geq 0$,
\begin{equation}
  \label{eq:pvpq}
  \begin{aligned}
    & \nu_{b+}^k - \nu_{b-}^k = v_b^k - v_b^0 \; ,\\
    & 0 \leq \nu_{b-}^k \perp \overline{q}_g - q_g^k \geq 0 \; , \\
    & 0 \leq \nu_{b+}^k \perp q_g^k - \underline{q}_g \geq 0 \; .
  \end{aligned}
\end{equation}
The corrective SCOPF modeling the recourse with \eqref{eq:droopmpec}
and \eqref{eq:pvpq} writes as a nonlinear program with complementarity
constraints. We note that the MPCC formulation has been adopted throughout the GO
competition, and is described at length in \cite{aravena2023recent,curtis2023decomposition,gholami2023admm}.

\subsection{AC-SCOPF}
We adapt the formulation of the AC-SCOPF used in \cite{frank2016introduction} to include the complementarity
constraints modeling the AGC~\eqref{eq:droopmpec} and the PV/QP switches~\eqref{eq:pvpq}.
We note $u_0 = (p_g^0, v_{b,PV}^0)$ and $x_0 = (\theta_b^0, v_{b,PQ}^0, q_g^0)$ the control and the state in the base case,
and $u_k = (p_g^k, v_{b,PV}^k)$ and $x^k = (\theta_b^k, v_{b,PQ}^k, q_g^k, \Delta^k, \pi_g^k, \nu^k_b)$
the control and the state in the contingency $k$. We define the AC-SCOPF problem as:
\begin{equation}
  \label{eq:scopf}
  \begin{aligned}
    \min_{x, u} \; & f(x_0, u_0) \\
    \text{s.t.} \quad & g_0(x_0, u_0) = 0 \; , \;  h_0(x_0, u_0) \leq 0 \; , \\
                      & \forall k \in \{1,\cdots,K\} :\\
                            & g_k(u_0, x_k, u_k) = 0 \;, \;  h_k(x_k, u_k) \leq 0 \; , \\
                            & 0 \leq G(x_k, u_k) \perp H(x_k, u_k) \geq 0 \; ,
  \end{aligned}
\end{equation}
where $u = (u_0, u_1, \cdots, u_K)$, $x = (x_0, x_1, \cdots, x_K)$,
$f(\cdot)$ is the objective in the base case scenario,
$g_k(\cdot)$ encodes the power-flow constraints and the two linear constraints
in \eqref{eq:droopmpec} and \eqref{eq:pvpq}, $h_k(\cdot)$ the operational
constraints (line-flow and operational bounds) and $G(\cdot)$ and
$H(\cdot)$ are two functions encoding the
complementarity parts in \eqref{eq:droopmpec} and \eqref{eq:pvpq}.

For a given base-case control $u_0$, the contingency $k$ is feasible if there is a solution
$(x_k, u_k)$ to the nonlinear system with complementarity constraints:
\begin{equation}
  \label{eq:screening}
  \left\{
  \begin{aligned}
    &  g_k(u_0, x_k, u_k) = 0 \; ,\;  h_k(x_k, u_k) \leq 0 \; , \\
    & 0 \leq G(x_k, u_k) \perp H(x_k, u_k) \geq 0 \; .
  \end{aligned}
  \right.
\end{equation}
The goal of \eqref{eq:scopf} is to find a base-case control $u_0$ such
that \eqref{eq:screening} is feasible for all the contingency $k=1,\cdots,K$.

\section{Mathematical programs with complementarity constraints}
We detail the formulation of the SCOPF~\eqref{eq:scopf} as a mathematical
program with complementarity constraints (MPCC) in \S\ref{sec:mpcc:intro},
and discuss the associated first-order complementarity conditions in \S\ref{sec:mpcc:stationary}.
Classical solution methods are discussed in \S\ref{sec:mpcc:solution}.

\subsection{MPCC in vertical formulation}
\label{sec:mpcc:intro}
Upon introducing slack variables to reformulate the inequality constraints
and the nonlinear complementarity constraints in \eqref{eq:scopf}, we obtain
the equivalent MPCC in vertical form:
\begin{equation}
  \label{eq:mpcc}
  \begin{aligned}
    \min_{w \in \mathbb{R}^n} \; & \phi(w) \quad \text{s.t.} \quad \left\{
      \begin{aligned}
        & c(w) = 0 \; , w_0 \geq 0 \; , \\
        & 0 \leq w_1 \perp w_2 \geq 0 \; ,
      \end{aligned}
    \right.
  \end{aligned}
\end{equation}
with
$w \in \mathbb{R}^n$ the variable aggregating the controls, states and slack variables
and $c: \mathbb{R}^n \to \mathbb{R}^m$ a function aggregating all the equality and inequality constraints
in \eqref{eq:scopf}.
We partition the decision variable as $w = (w_0, w_1, w_2) \in \mathbb{R}^{n-2p} \times
\mathbb{R}^p \times \mathbb{R}^p$ to isolate the variables contributing to the $p$ complementarity constraints.
For multipliers $(\lambda, \xi, \mu_1, \mu_2) \in \mathbb{R}^{m} \times \mathbb{R}^{n-2p} \times \mathbb{R}^p \times \mathbb{R}^p$, we define
the MPCC Lagrangian as
\begin{multline}
  \mathcal{L}^{\mathrm{MPCC}}(w, \lambda, \xi, \mu) = \phi(w) + \lambda^\top c(w) \\ - \xi^\top w_0 - \mu_1^\top w_1 - \mu_2^\top w_2 \; .
\end{multline}
The MPCC~\eqref{eq:mpcc} is equivalent to the nonlinear program:
\begin{equation}
  \label{eq:mpccnlp}
  \begin{aligned}
    \min_{w \in \mathbb{R}^n} \; & \phi(w) \quad \text{s.t.} \quad
    \left\{
    \begin{aligned}
      & c(w) = 0 \; , \; w_0 \geq 0 \; , \\
      & (w_1, w_2) \geq 0 \; , \; W_1 W_2 e \leq 0 \; ,
    \end{aligned}
    \right.
  \end{aligned}
\end{equation}
where we have noted $W_1 = \text{diag}(w_1)$ and $W_2 = \text{diag}(w_2)$ and $e \in \mathbb{R}^p$
a vector of ones.
The Lagrangian for \eqref{eq:mpccnlp} is
\begin{multline}
  \mathcal{L}(w, \lambda, \xi, \nu) = \phi(w) + \lambda^\top c(w) \\ - \xi^\top w_0 - \nu_1^\top w_1 - \nu_2^\top w_2 + \nu_0^\top W_1 W_2 e \; .
\end{multline}
The problem~\eqref{eq:mpccnlp} is degenerate, in the sense that
the relative interior of the feasible set is empty. As a consequence, \eqref{eq:mpccnlp}
does not satisfy the Mangasarian-Fromovitz constraint qualification (MFCQ), implying the multipliers $(\lambda, \nu)$ can be unbounded.
This can cause serious numerical issues if \eqref{eq:mpccnlp} is directly treated
by classical nonlinear programming solvers.

\subsection{First-order stationary conditions}
\label{sec:mpcc:stationary}
We note the feasible set of \eqref{eq:mpcc} as
\begin{equation*}
  \Omega = \{ w \in \mathbb{R}^n \; | \; c(w) = 0 \, , w_0 \geq 0 \; ,\; 0 \leq w_1 \perp w_2 \geq 0 \} \; .
\end{equation*}
For any feasible point $w \in \Omega$, we define the index sets
\begin{equation}
  \begin{aligned}
    \mathcal{I}_{+0}(w) = \{ i \in \{1, \cdots, p\} \; | \; w_{1,i} > 0 \, , \, w_{2,i} = 0 \} \;, \\
    \mathcal{I}_{0+}(w) = \{ i \in \{1, \cdots, p\} \; | \; w_{1,i} = 0 \, , \, w_{2,i} > 0 \} \;, \\
    \mathcal{I}_{00}(w) = \{ i \in \{1, \cdots, p\} \; | \; w_{1,i} = 0 \, , \, w_{2,i} = 0 \} \; .
  \end{aligned}
\end{equation}
The point $w \in \Omega$ satisfies \emph{strong stationarity}~\cite{scheel2000mathematical} if there exists
multipliers $(\lambda, \xi, \mu) \in \mathbb{R}^m \times \mathbb{R}^{2p}$ such that
\begin{equation}
  \label{eq:strongstationarity}
  \begin{aligned}
    & \nabla_w \mathcal{L}^{\mathrm{MPCC}}(w, \lambda, \xi, \mu) = 0 \; ,\\
    & c(w) = 0 \; , \; w_0 \geq 0 \; ,\\
    & w_{1,i} \geq 0 \; , \; \mu_{1,i} = 0 \;,\; w_{2, i} = 0 \;,\; \mu_{2,i} \in \mathbb{R} \;,\; \forall i \in \mathcal{I}_{+0}(w) \; , \\
    & w_{1,i} = 0 \; , \; \mu_{1,i} \in \mathbb{R} \;,\; w_{2, i} \geq 0 \;,\; \mu_{2,i} = 0 \;,\; \forall i \in \mathcal{I}_{0+}(w) \; , \\
    & w_{1,i} = 0 \; , \; \mu_{1,i} \geq 0 \;,\; w_{2, i} = 0 \;,\; \mu_{2,i} \geq 0 \;,\; \forall i \in \mathcal{I}_{00}(w) \; .
  \end{aligned}
\end{equation}
We note that the conditions
\eqref{eq:strongstationarity} are the KKT stationary solution of a relaxed nonlinear program~\cite{scheel2000mathematical} with no complementarity constraint.
If the relaxed nonlinear program satisfies Linear Independence Constraint Qualification (LICQ),
we say that MPCC-LICQ holds at $x \in \Omega$.
Our goal is to find a strong stationary solution for \eqref{eq:mpcc}.

\subsection{Solution methods}
\label{sec:mpcc:solution}
The solution of MPCC~\eqref{eq:mpcc} has been widely studied since the 2000s.
\emph{Direct methods} solves the problem~\eqref{eq:mpccnlp} directly using
a sequential quadratic programming (SQP) or an interior-point method (IPM)~\cite{fletcher2006local}.
\emph{Regularization methods} relax the degenerate terms in \eqref{eq:mpccnlp}
with a small parameter $\tau > 0$, such that the constraint $W_1 W_2 e \leq 0$
is reformulated as $W_1 W_2 e \leq \tau$. We obtain the so-called Scholtes
relaxation~\cite{scholtes2001convergence}: by driving the term $\tau$ to $0$ we recover the solution of the original MPCC~\eqref{eq:mpcc}.
\emph{Penalty-based methods} use a $\ell_1$-exact penalty to penalize the complementarity
constraints in the objective by adding a term $\tau w_1^\top w_2$, with a large-enough penalty $\tau >0$.
Assuming MPCC-LICQ, the $\ell_1$-exact penalty method converges to a strong-stationary solution $x \in \Omega$
\cite{ralph2004some}.
This is the method used in the solver Knitro~\cite{leyffer2006interior}.
In the next section, we interpret NCL as a mix of regularization and penalty-based methods.

\section{Augmented Lagrangian for MPCCs}
In this section, we present Algorithm NCL in \S\ref{sec:ncl:ncl} and
we analyze the structure of the Newton systems in \S\ref{sec:ncl:newton}.
We show that NCL can solve the problem~\eqref{eq:mpccnlp} on the GPU.

\subsection{Algorithm NCL}
\label{sec:ncl:ncl}
NCL is strictly equivalent to the classical augmented Lagrangian method,
but uses a nonlinearly constrained formulation in the subproblems.
When solving the MPCC~\eqref{eq:mpcc}, it has been proven in \cite{izmailov2012global} that the augmented Lagrangian
method converges to a strongly stationary solution $x \in \Omega$ if MPCC-LICQ holds at $x$
and the sequence of multiplier estimates $\{\nu_0^n\}_n$ generated by the algorithm
has a bounded subsequence (this prevents the algorithm to converge to a spurious solution).
NCL operates at two level: it updates the penalty $\rho^{(n)}$
and the two multiplier estimates $(\lambda^{(n)}, \nu_0^{(n)}) \in \mathbb{R}^m \times \mathbb{R}^p$ in the \emph{outer iterations}, whereas the \emph{inner iterations}
solves the constrained nonlinear subproblem:
\begin{equation}
  \label{eq:nclsubpb}
  \begin{aligned}
    \min_{w, r, t} \; & \; \mathcal{L}_\rho(w, r, t, \lambda^{(n)}, \nu_0^{(n)}) \\
    \text{s.t.} \quad & c(w) - r = 0 \; , \\
                      & W_1 W_2 e \leq t  \; , \; (w_0, w_1, w_2) \geq 0 \; ,
  \end{aligned}
\end{equation}
with
$\mathcal{L}_\rho(w, r, t, \lambda^{(n)}, \nu^{(n)}_0) =  \phi(w) + (\lambda^{(n)})^\top r  + (\nu_0^{(n)})^\top t +  \frac{\rho^{(n)}}{2}(\|r \|^2 + \|t\|^2)$
and $(r, t) \in \mathbb{R}^m \times \mathbb{R}^p$ two regularization variables.
The subproblem~\eqref{eq:nclsubpb} is always feasible.
We emphasize that NCL treats \eqref{eq:mpccnlp} as a generic nonlinear program, and does not
apply any special treatment to the complementarity constraints. The regularization
$t$ plays the role of the Scholtes relaxation parameter (see \S\ref{sec:mpcc:solution})
and is penalized in the objective by the Augmented Lagrangian penalty.

At a given outer iteration $n$, NCL solves the subproblem \eqref{eq:nclsubpb}
with an interior-point method (IPM).
Only the objective changes between two successive outer iterations, hence the IPM
can be efficiently warmstarted from the previous primal-dual solution.
Once the subproblem~\eqref{eq:nclsubpb} solved, we update the augmented Lagrangian
parameters $(\rho^{(n)}, \lambda^{(n)}, \nu_0^{(n)})$ using the solution $(x^{n+1}, r^{n+1}, t^{n+1})$
and proceed to the next outer iteration. We refer to \cite{montoison2025madncl} for a comprehensive
description of the algorithm.


If the problem is infeasible, the penalty $\rho^{(n)}$ diverges to $+\infty$
and NCL iterates are converging to a stationary solution of the \emph{feasibility problem}:
\begin{equation}
  \label{eq:feaspb}
  \begin{aligned}
    \min_{w, r, t} \; & \; \|r\|^2 + \|t\|^2 \\
    \text{s.t.} \quad & c(w) - r = 0 \\
                      & (w_0, w_1, w_2) \geq 0 \; , \; W_1 W_2 e \leq t \;.
  \end{aligned}
\end{equation}
On the contrary to the formulation used by the GO-competition~\cite{aravena2023recent}, we do not
have to penalize the constraint violation in the objective using a $\ell_1$ penalty, as the augmented
Lagrangian algorithm is doing that automatically.

\subsection{Newton systems}
\label{sec:ncl:newton}
The Karush-Kuhn-Tucker (KKT) stationary conditions for the subproblem~\eqref{eq:nclsubpb} are:
\begin{equation}
  \label{eq:kktncl}
  \begin{aligned}
    & \nabla \phi(w) + \nabla c(w)^\top \lambda = 0\\
    & \lambda^{(n)} + \rho^{(n)} r - \lambda = 0 \\
    & \nu_0^{(n)} + \rho^{(n)} t - \nu = 0 \\
    & c(w) - r = 0 \\
    & 0 \leq w_0 \perp \xi \geq 0 \\
    & 0 \leq t - W_1 W_2 e \perp \nu_0 \geq 0 \\
    & 0 \leq w_1 \perp \nu_1 \geq 0 \; , \; 0 \leq w_2 \perp \nu_2 \geq 0 \; .
  \end{aligned}
\end{equation}
Upon introducing a slack $s := t - W_1 W_2e$, a primal-dual interior-point method
solves the system~\eqref{eq:kktncl} for $(w, r, t, s, \lambda, \xi, \nu)$ using a Newton method. For a given
barrier parameter $\mu > 0$, the complementarity conditions in \eqref{eq:kktncl} are reformulated as
\begin{equation}
  W_0 \Xi = \mu e \; , \; S V_0 = \mu e \; , \; W_1 V_1 = \mu e \; , \; W_2 V_2 = \mu e \; ,
\end{equation}
with $\Xi = \text{diag}(\xi)$, $S = \text{diag}(s)$, $V_i = \text{diag}(\nu_i)$ for $i=0,1,2$.

Our implementation of NCL uses MadNLP to solve \eqref{eq:kktncl}. The barrier $\mu$ is updated
using the Fiacco-McCormick rule. The globalization is performed using a filter line-search
algorithm \cite{wachter2006implementation}. At each IPM iteration, the Newton system associated to NCL subproblem reduces to
\begin{equation}
  \label{eq:newtonsystem}
  \begin{bmatrix}
    A & \phantom{-}B^\top \\
    B & -C
  \end{bmatrix}
  \begin{bmatrix}
    \Delta w \\ \Delta y
  \end{bmatrix}
  =
  \begin{bmatrix}
    r_1 \\ r_2
  \end{bmatrix}
  \; ,
\end{equation}
with $\Delta y = (\Delta \lambda, \Delta \nu_0)$ the multipliers update
and $(r_1, r_2)$ an appropriate right-hand-side set by the IPM.
The matrix $A$ is defined for $H = \nabla^2_{ww} \mathcal{L}(w, r, t, \lambda,\xi,\nu)$:
\begin{equation*}
  \label{eq:blockA}
  A = \begin{bmatrix}
    H_{00} + W_0^{-1} \Xi & H_{01} & H_{02} \\
    H_{10} & H_{11} + W_1^{-1} V_1 & H_{12} + V_0 \\
    H_{20} & H_{21} + V_0 & H_{22} + W_2^{-1} V_2
  \end{bmatrix} \; ,
\end{equation*}
and with the Jacobian $J = \nabla c(w)^\top$ and regularization terms:
\begin{equation*}
  \label{eq:blockBC}
  B = \begin{bmatrix}
    J_0 & J_1 & J_2 \\
      & W_2 & W_1
  \end{bmatrix} \; , \;
  C = \begin{bmatrix}
    \rho^{-1} I &   \\
    & \rho^{-1} I + V_0^{-1} S
    \end{bmatrix} \; .
\end{equation*}
Once the system~\eqref{eq:newtonsystem} solved, the remaining
descent directions are recovered as
\begin{equation*}
  \begin{aligned}
    & \Delta \nu_i = W_i^{-1} (W_i V_ie - \mu e - V_i \Delta w_i)\;,   & i = 1,2 \; , \\
    & \Delta r = {\rho^{(n)}}^{-1} (\lambda + \rho^{(n)} r - \lambda^{(n)} - \Delta \lambda) \; ,\\
    & \Delta t = {\rho^{(n)}}^{-1} (\nu_0 + \rho^{(n)} t - \nu_0^{(n)} - \Delta \nu_0) \; .
  \end{aligned}
\end{equation*}

The Newton system~\eqref{eq:newtonsystem} reflects the structure of
the problem~\eqref{eq:nclsubpb}. Compared to a classical IPM method, the
$(2,2)$ block $C$ is non-zero: as it is well known, the augmented Lagrangian
term adds a natural dual regularization to the system. This regularization
accounts for NCL's robustness on degenerate instances.
Here, the MPCC structure leads to two potential degeneracies in \eqref{eq:newtonsystem}.
First, if strict complementarity does not hold
at the solution $w$ ($\mathcal{I}_{00}(w) \neq \emptyset$), the corresponding rows in $\begin{bmatrix}
  0 & W_2 & W_1 \end{bmatrix}$ converge to 0, implying $B$ is not full rank at the limit.
In that case, the non-zero block $C$ ensures that \eqref{eq:newtonsystem} remains non-singular.
Second, the bilinear terms $\begin{bmatrix} 0 & V_0 \\ V_0 & 0 \end{bmatrix}$
appearing in $A$ increases the indefinitess in \eqref{eq:newtonsystem},
but it does not impact the positive-definitess of the reduced Hessian
as we approach the solution~\cite[Section 3.4]{demiguel2005two}, keeping limited
the effect of this second degeneracy.

The linear system~\eqref{eq:newtonsystem} is sparse symmetric indefinite.
As a result, traditional optimization solvers often resort to HSL~\cite{duff1983multifrontal} to solve \eqref{eq:newtonsystem}.
These solvers use numerical pivoting for stability, a method known to be difficult to parallelize
and limiting the tractability of these solvers for large-scale SCOPF.
However, the non-zero block $C$ adds a regularization to the system, which allows
a solution of \eqref{eq:newtonsystem} with a pivoting-free factorization, e.g. using the signed Cholesky factorization~\cite{montoison2025madncl}.
Efficient parallel implementations of pivoting-free algorithms exist for GPUs,
opening the door for a fast factorization routine for \eqref{eq:newtonsystem}.

\section{Numerical results}
In this section, we test the performance of MadNCL on the AC-SCOPF problem~\eqref{eq:scopf}.
The implementation is discussed in \S\ref{sec:numerics:implementation}. We present
in \S\ref{sec:numerics:screening} how to perform the screening of the contingencies
with MadNCL, and in \S\ref{sec:numerics:scopf} the performance obtained by MadNCL when
solving AC-SCOPF instances.

\subsection{Implementation}
\label{sec:numerics:implementation}

MadNCL has been implemented in Julia 1.11, and uses the solver MadNLP to solve
the subproblems~\eqref{eq:nclsubpb}. As a consequence, MadNCL runs both on the CPU or on the GPU \cite{shin2024accelerating}.
The implementation is described in \cite{montoison2025madncl}
and is available at this URL: \url{github.com/MadNLP/MadNCL.jl}.
The CPU and GPU implementations of MadNCL are referred to as MadNCL-CPU and MadNCL-GPU, respectively.

In terms of running time, the two bottlenecks in MadNCL are (i) the evaluation of the derivatives
and (ii) the solution of the Newton system~\eqref{eq:newtonsystem}.
On the one hand, we have implemented the problem \eqref{eq:scopf} with ExaModels,
a modeler that exploits the repeated structures in the model to evaluate it in parallel.
We found that ExaModels is at least 10x faster than JuMP on the CPU~\cite{shin2024accelerating}.
Furthermore, ExaModels supports
the evaluation of the derivatives on the GPU, which gives us an additional speed-up factor.
On the other hand, the Newton systems~\eqref{eq:newtonsystem} are solved with HSL MA57 on the CPU.
On the GPU, we use the linear solver NVIDIA cuDSS, implementing a signed Cholesky factorization.

We compare MadNCL with the interior-point solver Knitro~\cite{waltz2006interior},
running on the CPU with the linear solver HSL MA57. Knitro evaluates the SCOPF model~\eqref{eq:scopf}
with JuMP. Knitro supports problems with complementarity constraints and implements the $\ell_1$-exact penalty
method described in \cite{leyffer2006interior}.
As a contender in the GO-competition, the support of complementarity constraints in Knitro has
been significantly improved in the recent years.

We use instances from the MATPOWER library~\cite{zimmerman2010matpower}.
In the following experiments the processor is an AMD EPYC 7430.
The GPU is a NVIDIA A30 with 24GB of device memory.
We provide the code to reproduce the benchmark in this repository: \url{github.com/frapac/pscc-scopf}.
In all our experiments, we set the convergence tolerance to {\tt tol=1e-6}.
In MadNCL, we set the minimum barrier parameter to $\mu_{min} = 10^{-7}$ inside MadNLP's IPM algorithm.

\subsection{Contingency screening}
\label{sec:numerics:screening}
Screening of contingencies is an important part in SCOPF's pre-processing~\cite{capitanescu2007contingency}.
The goal is to reduce the total number of contingencies $K$ in \eqref{eq:scopf} by identifying a set of representative critical contingencies.
Unfortunately, not all the contingencies are feasible: some
are \emph{structurally infeasible} (there is no base solution $u_0$ such that
\eqref{eq:screening} is feasible), and some others are conflicting (for two contingencies
$k, l$, there is no $u_0$ such that the two associated problems~\eqref{eq:screening} are feasible).
In that case, NCL falls back automatically to the feasibility problem~\eqref{eq:feaspb}
by increasing the penalty $\rho^{(n)}$ to infinity. As a result, the algorithm returns the
regularizations $(r, t)$ minimizing the local infeasibility (with guarantees in the convex case~\cite{chiche2016augmented}).

\begin{figure}[!ht]
\centering
\includegraphics[width=.4\textwidth]{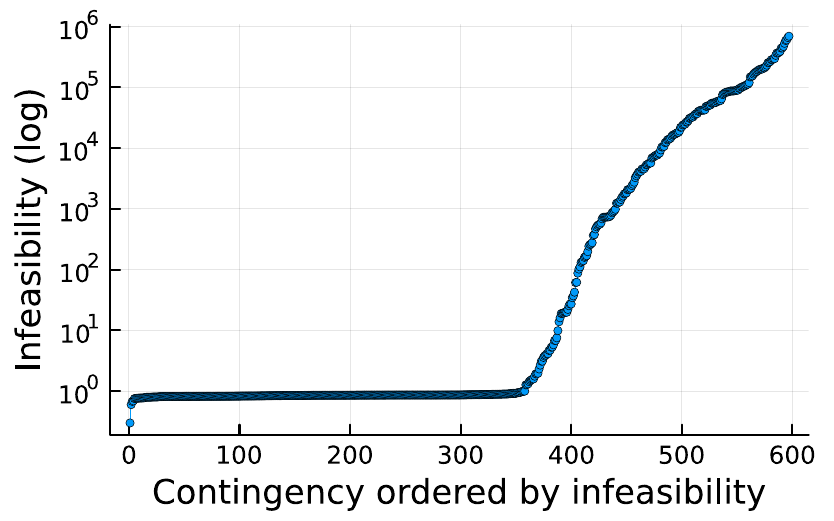}
\caption{Result of contingency screening for {\tt ACTIVSg500}.
The base case solution in \eqref{eq:screening} is obtained by solving
the classical OPF problem.}
\label{fig:contingency}
\end{figure}

We assess the performance of MadNCL when scanning the contingencies in a given network.
We compute a base case solution $u_0$ by solving the AC OPF problem, and solve the
system~\eqref{eq:screening} using MadNCL. MadNCL reformulates
the system~\eqref{eq:screening} as a feasibility problem~\eqref{eq:feaspb}.
The final objective value quantifies the contingency's infeasibility: the closer to 0,
the closer to feasibility. On the contrary, a large objective value indicates that the contingency is infeasible.
We compare the average time to screen one line contingency in Table~\ref{tab:timecontingency}.
Both Knitro and MadNCL solves the system~\eqref{eq:screening} and run on the CPU using the linear solver HSL MA57. Knitro handles explicitly the complementarity constraints, whereas MadNCL
looks at the NLP~\eqref{eq:mpccnlp}. We observe that in term
of computation time, MadNCL is at least 10 times faster than Knitro on instances
with more than 300 buses: MadNCL is able to detect
infeasible problems (a common occurrence in contingency screening) much faster than
the algorithm implemented in Knitro.
We display the result returned by the algorihtm in Figure~\ref{fig:contingency}: we scan all the line contingencies for
the instance {\tt ACTIVSg500} and we order them by level of infeasibility (as measured
by the final objective returned by MadNCL).

\begin{table}[!ht]
\centering
\caption{Average time to scan one contingency (in seconds) for different instances.}
\begin{tabular}{|r|rr|}
  \hline
  Case  & Knitro (s) & MadNCL-CPU (s) \\
  \hline
  118ieee & 0.5 & 0.01 \\
  ACTIVSg200 & 0.5 & 0.1 \\
  300ieee & 5.5 & 0.2 \\
  ACTIVSg500 & 5.4 &  0.3  \\
  1354pegase & 75.4  & 5.0 \\
  ACTIVSg2000 & 40.7 & 2.5  \\
  2869pegase & 238.4  & 14.1  \\
  \hline
\end{tabular}
\label{tab:timecontingency}
\end{table}

\subsection{Solution of large-scale AC-SCOPFs on GPU}
\label{sec:numerics:scopf}
We are now in position to answer the main research question in this paper: is
MadNCL effective at solving large-scale AC-SCOPFs on the GPU?
We set up the following experiment.
We compare Knitro and MadNCL-CPU (with HSL MA57) with MadNCL-GPU (using NVIDIA cuDSS).
Knitro uses the modeler JuMP, which supports passing the complementarity constraints explicitly to the solver.
On the contrary, MadNCL converts the JuMP model to ExaModels~\cite{shin2024accelerating} to benefit from faster evaluation
of derivatives (in particular on the GPU). MadNCL solves the NLP reformulation~\eqref{eq:mpccnlp}.
We select $K$ representative contingencies using the method described in \S\ref{sec:numerics:screening},
and ensure that all the contingencies are not structurally infeasible (so \eqref{eq:scopf} remains feasible).

\begin{figure}[!ht]
\centering
\includegraphics[width=.4\textwidth]{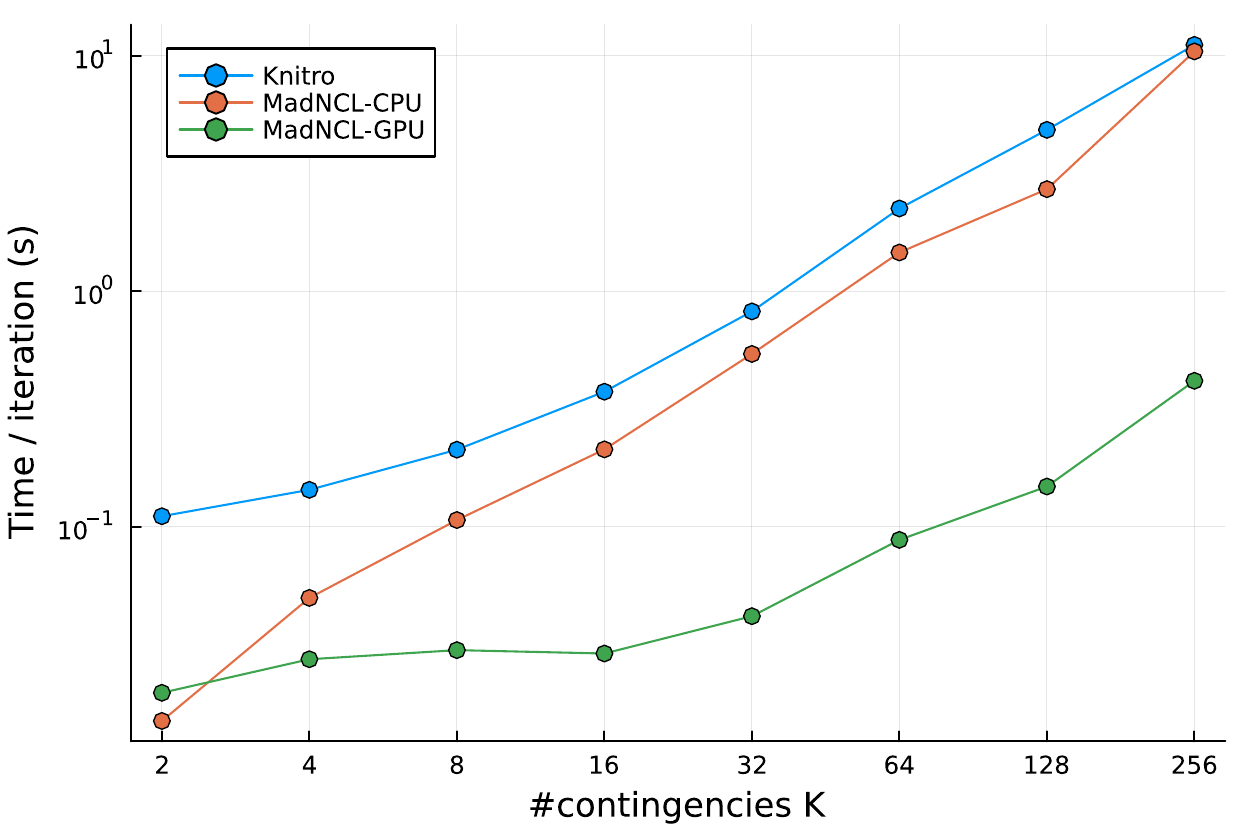}
\caption{Time per iteration (in seconds) against number of contingencies $K$
  when solving {\tt ACTIVSg500}. We compare Knitro (using JuMP)
against MadNCL (using ExaModels). We test MadNCL both on the CPU (using HSL MA57)
and on the GPU (using NVIDIA cuDSS).}
\label{fig:scalability}
\end{figure}

First, we consider the case {\tt ACTIVSg500} and increase the number of contingencies
$K$ from $2$ to $256$ in \eqref{eq:scopf}: with $K=256$, the SCOPF~\eqref{eq:scopf} has almost 1 million of variables.
We detail the results in Table~\ref{tab:scopfscalability}.
Interestingly, the number of iterations in MadNCL-CPU and MadNCL-GPU is different: as we enter in the large-scale
regime, HSL MA57 and NVIDIA cuDSS are returning slightly different solutions, resulting in small discrepancies in the algorithm.
This is amplified by the parallel nature of NVIDIA cuDSS: consecutive runs in NVIDIA cuDSS can be non-deterministic for large $K$.
The objective value gives an insight on the quality of the solution: as the problem \eqref{eq:scopf} is non-convex, different algorithms can converge
to different solutions. However, we observe here that MadNCL-CPU and MadNCL-GPU both converges to the same solution as Knitro
We show in Figure~\ref{fig:scalability} the time per iteration (in seconds) against
the number of contingencies for Knitro, MadNCL-CPU and MadNCL-GPU.
In term of raw performance, MadNCL-CPU is slightly faster than Knitro: MadNCL
here benefits from faster evaluations of the derivatives with ExaModels (compared to JuMP).
MadNCL-GPU decreases the time per iteration further by at least an order of magnitude:
for $K \geq 64$, MadNCL-GPU becomes at least 20x faster than Knitro and 6x faster than MadNCL-CPU. The linear solver cuDSS is
here very effective at solving the linear system \eqref{eq:newtonsystem}, explaining the faster solution time on the GPU.

\begin{table*}[!ht]
\centering
\caption{Performance of Knitro and MadNCL as we increase the number
of contingencies $K$ for {\tt ACTIVSg500}. We display the total number of IPM iterations
and the time to solution (in seconds). The symbol "-" indicates that the solver has failed
to find a solution.}
\begin{tabular}{|rrr|rrr|rrr|rrr|}
\hline
\multicolumn{ 3}{|c|}{}  & \multicolumn{ 3}{c|}{\bf Knitro} & \multicolumn{ 3}{c|}{\bf MadNCL-CPU} & \multicolumn{ 3}{c|}{\bf MadNCL-GPU}  \\
\hline
K & nvar & ncon & Iter & Obj. & Time (s)& Iter & Obj. & Time (s)& Iter & Obj. & Time (s) \\
\hline
4 & 18400 & 24251 & 355 & 7.28 & 51.01 & 238 & 7.28 & 7.31 & 240 & 7.28 & 4.45\\
8 & 33300 & 43919 & 418 & 7.28 & 88.94 & 435 & 7.28 & 29.23 & 290 & 7.28 & 7.77\\
16 & 63100 & 83255 & 114 & 7.28 & 42.76 & 214 & 7.28 & 25.71 & 261 & 7.28 & 6.65\\
32 & 122700 & 161927 & 345 & 7.28 & 283.68 & 587 & 7.28 & 166.20 & 568 & 7.28 & 23.64\\
64 & 241900 & 319271 & 960 & 7.28 & 2159.59 & 528 & 7.28 & 273.08 & 453 & 7.28 & 27.96\\
128 & 480300 & 633959 & - & 7.29 & 4852.33 & 415 & 7.29 & 421.09 & 265 & 7.29 & 46.40\\
256 & 957100 & 1263335 & - & 7.30 & 11136.08 & 493 & 7.30 & 1120.16 & 609 & 7.30 & 170.75\\
\hline
\end{tabular}
\label{tab:scopfscalability}
\end{table*}

Second, we compare Knitro and MadNCL-GPU on different instances, with a varying number of contingencies $K$.
The results are displayed in Table~\ref{tab:benchmark}.
We observe that MadNCL-GPU is consistently faster than Knitro. As before in Table~\ref{tab:scopfscalability},
we do not have any guarantee that Knitro converges to the same solution as MadNCL: they
report a different solution for {\tt 1354pegase} and {\tt 2869pegase}.
We notice that MadNCL converges in significantly more iterations
for {\tt ACTIVSg2000} and {\tt 2869pegase}. This points to MadNCL's main limitation:
as we increase the problem's dimension, we increase the degeneracy in the indefinite linear system~\eqref{eq:newtonsystem}.
As a result, the filter line-search algorithm used internally to solve the subproblem~\eqref{eq:nclsubpb} has to
perform numerous primal-dual regularizations during the inertia correction, significantly impairing MadNCL's convergence (we have observed that the ill-conditioning arises mostly from the reformulation
of \eqref{eq:pvpq}, when the reactive power lower-bound $\underline{q}_g$ is close to the upper-bound $\overline{q}_g$).
For that reason, we have observed that MadNCL does not converge consistently on larger instances.
We plan to address this issue in future work.

\begin{table*}[!ht]
\centering
\caption{Performance of Knitro and MadNCL-GPU on different instances. We
  display the total number of IPM iterations and the time to solution (in
  seconds). The symbol "-" indicates that the solver has failed to find a solution.
}
\begin{tabular}{|rr|rr|rrr|rrr|}
\hline
\multicolumn{ 4}{|c|}{}  & \multicolumn{ 3}{c|}{\bf Knitro} & \multicolumn{ 3}{c|}{\bf MadNCL-GPU}  \\
\hline
Name & K & nvar & ncon
& Iter & Obj. & Time (s)& Iter & Obj. & Time (s) \\
\hline
ACTIVSg200 & 10 & 17546 & 22841 & 141 & 2.76 & 14.84 & 59 & 2.76 & 3.58\\
ACTIVSg200 & 50 & 81906 & 106721 & 147 & 2.76 & 106.99 & 167 & 2.76 & 5.01\\
ACTIVSg200 & 100 & 162356 & 211571 & 81 & 2.76 & 97.93 & 244 & 2.76 & 13.69\\
\hline
ACTIVSg500 & 10 & 40750 & 53753 & 290 & 7.47 & 82.27 & 130 & 7.47 & 3.51\\
ACTIVSg500 & 50 & 189750 & 250433 & 533 & 7.83 & 871.45 & 366 & 7.83 & 34.62\\
ACTIVSg500 & 100 & 376000 & 496283 & 294 & 7.83 & 935.19 & 393 & 7.83 & 39.33\\
\hline
1354pegase & 8 & 109056 & 144327 & 43 & 7.41 & 53.82 & 124 & 7.42 & 7.66\\
1354pegase & 16 & 206920 & 273999 & 30 & 7.41 & 86.90 & 111 & 7.41 & 10.60\\
1354pegase & 32 & 402648 & 533343 & 116 & 7.41 & 656.68 & 410 & 7.42 & 83.86\\
\hline
ACTIVSg2000 & 8 & 173024 & 229853 & 160 & 122.89 & 1442.62 & 989 & 122.89 & 129.80\\
ACTIVSg2000 & 16 & 328360 & 436469 & 141 & 122.89 & 3001.10 & - & 122.90 & 220.79\\
\hline
2869pegase & 8 & 242102 & 323479 & 65 & 13.40 & 308.77 & 334 & 13.42 & 41.08\\
2869pegase & 16 & 459118 & 613727 & 68 & 13.40 & 775.95 & 952 & 13.42 & 225.56\\
\hline
\end{tabular}
\label{tab:benchmark}
\end{table*}

\section{Conclusion}

In this article, we have studied the solution of large-scale AC-SCOPFs
using MadNCL, a GPU-accelerated solver. As illustrated by Figure~\ref{fig:scalability},
MadNCL is significantly faster on the GPU than a similar solver running on the CPU,
with a time per iteration decreased by a factor of almost 40 on the largest instances.
As a result, MadNCL can find in less than 3 minutes a strongly stationary solution for AC-SCOPFs with almost a million of decision variables, allowing to tackle instances with thousands of buses and
hundred of contingencies. In future work, we plan to move that threshold further:
this would require addressing the ill-conditioning appearing inside the IPM's Newton systems,
caused by the millions of complementarity constraints found in AC-SCOPF with
more than 10,000 buses and 100 contingencies.

\end{document}